# GOODNESS-OF-FIT PROBLEM FOR ERRORS IN NONPARAMETRIC REGRESSION: DISTRIBUTION FREE APPROACH

By Estate V. Khmaladze and Hira L. Koul[1]

*Victoria University of Wellington and Michigan State University*

This paper discusses asymptotically distribution free tests for the classical goodness-of-fit hypothesis of an error distribution in nonparametric regression models. These tests are based on the same martingale transform of the residual empirical process as used in the one sample location model. This transformation eliminates extra randomization due to covariates but not due the errors, which is intrinsically present in the estimators of the regression function. Thus, tests based on the transformed process have, generally, better power. The results of this paper are applicable as soon as asymptotic uniform linearity of nonparametric residual empirical process is available. In particular they are applicable under the conditions stipulated in recent papers of Akritas and Van Keilegom and Müller, Schick and Wefelmeyer.

**1. Introduction.** Consider a sequence of i.i.d. pairs of random variables $\{(X_i, Y_i)_{i=1}^n\}$ where $X_i$ are $d$-dimensional covariates and $Y_i$ are the one-dimensional responses. Suppose $Y_i$ has regression in mean on $X_i$, that is, there is a regression function $m(\cdot)$ and a sequence of i.i.d. zero mean innovations $\{e_i, 1 \le i \le n\}$, independent of $\{X_i\}$, such that

$$Y_i = m(X_i) + e_i, \qquad i = 1, \ldots, n.$$

This regression function, as in most applications, is generally unknown and we do not make assumptions about its possible parametric form, so that we need to use a nonparametric estimator $\hat{m}_n(\cdot)$ based on $\{(X_i, Y_i)_{i=1}^n\}$.

The problem of interest here is to test the hypothesis that the common distribution function (d.f.) of $e_i$ is a given $F$. Since $m(\cdot)$ is unknown we can only use residuals

$$\hat{e}_i = Y_i - \hat{m}_n(X_i), \qquad i = 1, \ldots, n,$$

Received February 2008; revised September 2008.

[1]Supported in part by the NSF Grant DMS-07-04130.

which, obviously, are not i.i.d. anymore. Let $F_n$ and $\hat{F}_n$ denote the empirical d.f. of the errors $e_i, 1 \le i \le n$, and the residuals $\hat{e}_i, 1 \le i \le n$, respectively, and let

$$v_n(x) = \sqrt{n}[F_n(x) - F(x)], \qquad \hat{v}_n(x) = \sqrt{n}[\hat{F}_n(x) - F(x)], \qquad x \in \mathbb{R},$$

denote empirical and "estimated" empirical processes.

Akritas and Van Keilegom (2001) and Müller, Schick and Wefelmeyer (2007) established, under the null hypothesis and some assumptions and when $d = 1$, the following uniform asymptotic expansion of $\hat{v}_n$:

$$(1.1) \qquad \hat{v}_n(x) = v_n(x) - f(x)R_n + \xi_n(x), \qquad \sup_x |\xi_n(x)| = o_p(1),$$

where

$$(1.2) \qquad\qquad\qquad\qquad R_n = O_p(1).$$

Basically, the term $R_n$ is made up by the sum

$$R_n = n^{-1/2} \sum_{i=1}^{n} [\hat{m}_n(X_i) - m(X_i)],$$

but using special form of the estimator $\hat{m}_n$, Müller, Schick and Wefelmeyer obtained especially simple form for it:

$$(1.3) \qquad\qquad\qquad\qquad R_n = n^{-1/2} \sum_{i=1}^{n} e_i.$$

Müller, Schick and Wefelmeyer (2009) provides a set of sufficient conditions under which (1.1)–(1.3) continue to hold for the case $d > 1$.

In the case of parametric regression where the regression function is of the parametric form, $m(\cdot) = m(\cdot, \theta)$, and the unknown parameter $\theta$ is replaced by its estimator $\hat{\theta}_n$, similar asymptotic expansion have been established in Loynes (1980), Koul (2002) and Khmaladze and Koul (2004). However, the nonparametric case is more complex and it is remarkable that the asymptotic expansions (1.1) and (1.2) are still true.

The above expansion leads to the central limit theorem for the process $\hat{v}_n$, and, hence, produces the null limit distribution for test statistics based on this process. However, the same expansion makes it clear that the statistical inference based on $\hat{v}_n$ is inconvenient in practice and even infeasible; not only does the limit distribution of $\hat{v}_n$ after time transformation $t = F(x)$ still depend on the hypothetical d.f. $F$, but it depends also on the estimator $\hat{m}_n$ (and, in general, on the regression function $m$ itself), that is, it is different for different estimators. Since goodness-of-fit statistics are essentially nonlinear functionals of the underlying process with difficult to calculate limit distributions, it is practically inconvenient to be obliged to do substantial computational work to evaluate their null distributions every time we test the hypothesis. Note, in particular, that if we try to use some kind of bootstrap simulations, we would have to compute the nonparametric estimator $\hat{m}_n$ for every simulated subsample, which makes it especially time consuming.

Starting with asymptotic expansion (1.1) of Akritas and Van Keilegom and Müller, Schick and Wefelmeyer, our goal is to show that the above-mentioned complications can be avoided in the way, which is technically surprisingly simple. Namely, we present the transformed process $w_n$, which, after time transformation $t = F(x)$, converges in distribution to a standard Brownian motion, for any estimator $\hat{m}_n$ for which (1.1) is valid. One would expect that this is done at the cost of some power. We shall see, however, somewhat unexpectedly, that tests based on this transformed process $w_n$ should, typically, have better power than those based on $\hat{v}_n$. Perhaps it is worth emphasizing that to achieve this goal we actually need only the smallness of the remainder process $\xi_n$ and not asymptotic boundedness (1.2) in the expansion (1.1).

We end this section by mentioning some recent applications of martingale transform, in different types of regression problems, by Koenker and Xie (2002, 2006), Bai (2003), Delgado, Hidalgo and Velasco (2005) and Koul and Yi (2006).

**2. Transformed process.** Suppose the d.f. $F$ has an absolutely continuous density $f$ with a.e. derivative $\dot{f}$ and finite Fisher information for location. Let $\psi_f = -\dot{f}/f$ denote the score function for location family $F(\cdot - \theta), \theta \in \mathbb{R}$ at $\theta = 0$—we can assume that $\theta = 0$ without loss of generality. Then,

$$(2.1) \qquad \int \psi_f^2(x)\,dF(x) < \infty.$$

Consider augmented score function

$$h(x) = \begin{pmatrix} 1 \\ \psi_f(x) \end{pmatrix}$$

and augmented incomplete information matrix

$$\Gamma_{F(x)} = \int_x^\infty h(x)h^T(x)\,dF(x) = \begin{pmatrix} 1 - F(x) & f(x) \\ f(x) & \sigma_f^2(x) \end{pmatrix}, \qquad x \in \mathbb{R},$$

with $\sigma_f^2(x) = \int_x^\infty \psi_f^2(y)\,dF(y)$.

For a signed measure $\nu$ for which the following integral is well defined, let

$$K(x, \nu) = \int_{-\infty}^x h^T(y)\Gamma_{F(y)}^{-1} \int_y^\infty h(z)\,d\nu(z)\,dF(y), \qquad x \in \mathbb{R}.$$

Occasionally, $\nu$ will be a vector of signed measures in which case $K$ will be a vector also.

Our transformed process $w_n$ is defined as

$$(2.2) \qquad w_n(x) = \sqrt{n}[\hat{F}_n(x) - K(x, \hat{F}_n)], \qquad x \in \mathbb{R}.$$

We shall show that $w_n$ converges in distribution to the Brownian motion $w$ in time $F$, that is, $w_n(F^{-1})$ converges weakly to standard Brownian motion on the interval $[0, 1]$, where $F^{-1}(u) = \inf\{x; F(x) \geq u\}$, $0 \leq u \leq 1$.

To begin with observe that the process $w_n$ can be rewritten as

$$(2.3) \qquad w_n(x) = \hat{v}_n(x) - K(x, \hat{v}_n).$$

Indeed, $F(x)$ is the first coordinate of the vector-function $H(x) = \int_{-\infty}^{x} h \, dF = (F(x), -f(x))^T$, and we will see that

$$(2.4) \qquad H^T(x) - K(x, H^T) = 0 \qquad \forall x \in \mathbb{R}.$$

Subtracting this identity from (2.2) yields (2.3). Using asymptotic expansion (1.1) we can rewrite

$$(2.5) \quad w_n(x) = v_n(x) - K(x, v_n) + \eta_n(x), \qquad \eta_n(x) = \xi_n(x) - K(x, \xi_n),$$

where one expects $\eta_n$ to be "small" (see Section 4), and the main part on the right not to contain the term $f(F^{-1}(t))R_n$ of that expansion. This is true again because of (2.4) and because the second coordinate of $H(x)$ is $-f(x)$.

The transformation $w_n$ is very similar to the one studied in Khmaladze and Koul (2004) where regression function is assumed to be parametric. However, asymptotic behavior of the empirical distribution function $\hat{F}_n$ here is more complicated. As a result, we have to prove the smallness of the "residual process" $\eta_n$ in (2.5) differently (see Section 4). Here we demonstrate that although, in this transformation, singularity at $t = 1$ exists, the process $w_n(F^{-1})$ converges to its weak limit on the closed interval $[0, 1]$—see Theorem 4.1(ii). Besides, we explicitly consider the case of possibly degenerate matrix $\Gamma_{F(x)}$ and show that $w_n$ is still well defined—see Lemma 2.1.

If $\Gamma_{F(x)}$ is of the full rank for all $x \in \mathbb{R}$, then (2.4) is obvious. For most d.f.'s $F$, the matrix $\Gamma_{F(x)}$ indeed is not degenerate, that is, the coordinates $1$ and $\psi_f$ of $h$ are linearly independent functions on tail set $\{x > x_0\}$ for every $x_0 \in \mathbb{R}$. However, if (and only if) for $x$ greater than some $x_0$, the density $f$ has the form $f(x) = \alpha e^{-\alpha x}, \alpha > 0$, the function $\psi_f(x)$ equals the constant $\alpha$ so that $1$ and $\psi_f(x)$ become linearly dependent for $x > x_0$. As this can indeed be the case in applications, for example, for the double exponential distribution, it is useful to show that (2.4) is still correct and the transformation (2.3) still can be used.

The lemma below shows, that although in this case $\Gamma_{F(x)}^{-1}$ cannot be uniquely defined, the function $h^T(x) \Gamma_{F(x)}^{-1} \int_x^\infty h(y) \, d\mu(y)$ with $\mu = v_n$ or $\mu = \hat{v}_n$, is well defined. Here it is more transparent and simple to use also time transformation $t = F(x)$. Accordingly, let $u_n(t) = v_n(F^{-1}(t))$, $\hat{u}_n(t) = \hat{v}_n(F^{-1}(t))$, $\gamma(t) = h(F^{-1}(t))$, and $\Gamma_t = \int_t^1 \gamma(s)\gamma(s)^T \, ds$, $0 \leq t \leq 1$.

LEMMA 2.1. *Suppose, for some $x_0$, such that $0 < F(x_0) < 1$, the matrix $\Gamma_{F(x)}$, for $x > x_0$ degenerates to the form*

$$(2.6) \qquad \Gamma_{F(x)} = (1 - F(x)) \begin{pmatrix} 1 & \alpha \\ \alpha & \alpha^2 \end{pmatrix} \qquad \forall x > x_0, \text{ some } \alpha > 0.$$

*Then, the equalities (2.4) and, hence, (2.3) are still valid. Besides,*

$$h^T(x) \Gamma_{F(x)}^{-1} \int_x^\infty h(y) \, dv_n(y) = -\frac{v_n(x)}{1 - F(x)} \qquad \forall x \in \mathbb{R},$$

*or*

$$\gamma^T(t) \Gamma_t^{-1} \int_t^1 \gamma(s) \, du_n(s) = -\frac{u_n(t)}{1 - t} \qquad \forall 0 \le t < 1.$$

*A similar fact holds with $v_n(u_n)$ replaced by $\hat{v}_n(\hat{u}_n)$.*

REMARK 2.1. The argument that follows is an adaptation and simplification of a general treatment of the case of degenerate matrices $\Gamma_{F(x)}$, given in Nikabadze (1987) and Tsigroshvili (1998).

PROOF OF LEMMA 2.1. Let $\gamma(t) = (1, \alpha)^T, t = F(x)$. The image and kernel of the linear operator in $\mathbb{R}^2$ of $\Gamma_t$, respectively, are

$$\mathcal{I}(\Gamma_t) = \{b : b = \Gamma_t a \text{ for some } a \in \mathbb{R}^2\}$$
$$= \{b : b = \beta(1 - t)(1, \alpha)^T, \beta \in \mathbb{R}\};$$
$$\mathcal{K}(\Gamma_t) = \{a : \Gamma_t a = 0\} = \{a : a = c(-\alpha, 1)^T, c \in \mathbb{R}\}.$$

Moreover, both $\int_t^1 \gamma \, du_n$ and $H(F^{-1}(t))$ are in $\mathcal{I}(\Gamma_t)$ and if $b \in \mathcal{I}(\Gamma_t)$ then $\Gamma_t b = (1 - t)(1 + \alpha^2) b$. Then $\Gamma_t^{-1}$ is *any* (matrix of) linear operator on $\mathcal{I}(\Gamma_t)$ such that

$$\Gamma_t^{-1} b = \frac{1}{(1 - t)(1 + \alpha^2)} b + a, \qquad a \in \mathcal{K}(\Gamma_t).$$

But $\gamma(t) = (1, \alpha)^T$ is orthogonal to an $a \in \mathcal{K}(\Gamma_t)$ and therefore

$$(2.7) \qquad \gamma^T(t) \Gamma_t^{-1} b = \frac{1}{(1 - t)(1 + \alpha^2)} \gamma^T(t) b$$

does not depend on the choice of $a \in \mathcal{K}(\Gamma_t)$ and, hence, is defined uniquely. For $b = \int_t^1 \gamma(s) \, du_n(s)$ this gives the equality in the lemma. Besides, for any $b \in \mathcal{I}(\Gamma_t)$, $a \in \mathcal{K}(\Gamma_t)$,

$$\gamma^T(t) \Gamma_t^{-1} \Gamma_t(b + a) = \gamma^T(t) \Gamma_t^{-1} \Gamma_t b = \gamma^T(t) b = \gamma^T(t)(b + a),$$

which gives (2.4). The rest of the claim is obvious. □

Now consider the leading term of (2.5) in time $t = F(x)$. It is useful to consider its function parametric version, defined as

$$(2.8) \qquad b_n(\varphi) = u_n(\varphi) - K_n(\varphi), \qquad \varphi \in L_2[0,1],$$

where $u_n(\varphi) = \int_0^1 \varphi(s)\, du_n(s)$, and

$$K_n(\varphi) = K(\varphi, u_n) = \int_0^1 \varphi(t) \gamma^T(t) \Gamma_t^{-1} \int_t^1 \gamma(s)\, du_n(s)\, dt.$$

With slight abuse of notation, denote $b_n(\varphi)$ when $\varphi(\cdot) = I(\cdot \leq t)$ by

$$(2.9) \qquad b_n(t) = u_n(t) - \int_0^t \gamma^T(u) \Gamma_u^{-1} \int_u^1 \gamma(s)\, du_n(s)\, du.$$

Conditions for weak convergence of $u_n$ are well known: if $\Phi \subset L_2[0,1]$ is a class of functions, such that the sequence $u_n(\varphi), n \geq 1$, is uniformly in $n$ equicontinuous on $\Phi$, then $u_n \to_d u$ in $l_\infty(\Phi)$ where $u$ is standard Brownian bridge, see, for example, van der Vaart and Wellner (1996). The conditions for the weak convergence of $K_n$ to great extent must be simpler, because, unlike $u_n$, $K_n$ is continuous linear functional in $\varphi$ on the whole of $L_2[0,1]$, however, not uniformly in $n$. We will see, Proposition 2.1 below, that although, for every $\varepsilon > 0$, the provisional limit in distribution of $K_n(\varphi)$, namely,

$$K(\varphi) = K(\varphi, u) = \int_0^1 \varphi(t) \gamma^T(t) \Gamma_t^{-1} \int_t^1 \gamma(s)\, du(s)\, dt$$

is continuous on $L_{2,\varepsilon}$, the class of functions in $L_2[0,1]$ which are equal 0 on the interval $(1 - \varepsilon, 1]$, it is not continuous on $L_2[0,1]$. Therefore it is unavoidable to use some condition on $\varphi$ at $t = 1$. Condition (2.10) below still allows $\varphi(t) \to \infty$ as $t \to 1$ (see examples below).

THEOREM 2.1. (i) *Let $L_{2,\varepsilon} \subset L_2[0,1]$ be the subspace of all square integrable functions which are equal to 0 on the interval $(1 - \varepsilon, 1]$. Then, $K_n \to_d K$, on $L_{2,\varepsilon}$, for any $0 < \varepsilon < 1$.*

(ii) *Let, for an arbitrary small but fixed $\varepsilon > 0$, $C < \infty$, and $\alpha < 1/2$, $\Phi_\varepsilon \subset L_2[0,1]$ be a class of all square integrable functions satisfying the following right tail condition:*

$$(2.10) \qquad |\varphi(s)| \leq C[\gamma^T(s) \Gamma_s^{-1} \gamma(s)]^{-1/2} (1-s)^{-1/2-\alpha} \qquad \forall s > 1 - \varepsilon.$$

*Then, $K_n \to_d K$, on $\Phi_\varepsilon$.*

PROOF. (i) The integral $\int_t^1 \gamma\, du_n$ as process in $t$, obviously, converges in distribution to the Gaussian process $\int_t^1 \gamma\, du$. Therefore, all finite-dimensional distributions of $\gamma^T(t) \Gamma_t^{-1} \int_t^1 \gamma\, du_n$, for $t < 1$, converge to corresponding finite-dimensional distributions of the Gaussian process $\gamma^T(t) \Gamma_t^{-1} \int_t^1 \gamma\, du$. Hence, for any fixed $\varphi \in L_{2,\varepsilon}$, distribution of $K_n(\varphi)$ converges to that of $K(\varphi)$. So, we only need to show tightness, or, equivalently, equicontinuity of $K_n(\varphi)$ in $\varphi$. We have

$$|K_n(\varphi)| \leq \int_0^1 |\varphi(t)| \left| \gamma^T(t) \Gamma_t^{-1} \int_t^1 \gamma(s) \, du_n(s) \right| dt$$

$$\leq \sup_{t \leq 1-\varepsilon} \left| \gamma^T(t) \Gamma_t^{-1} \int_t^1 \gamma(s) \, du_n(s) \right| \int_0^{1-\varepsilon} |\varphi(t)| \, dt,$$

while

$$\sup_{t \leq 1-\varepsilon} \left| \gamma^T(t) \Gamma_t^{-1} \int_t^1 \gamma(s) \, du_n(s) \right| \to_d \sup_{t \leq 1-\varepsilon} \left| \gamma^T(t) \Gamma_t^{-1} \int_t^1 \gamma(s) \, du(s) \right| = O_p(1).$$

This proves that $K_n(\varphi)$ is equicontinuous in $\varphi \in L_{2,\varepsilon}$ and (i) follows.

(ii) To prove (ii), what we need is to show the equicontinuity of $K_n(\varphi)$ on $\Phi_\varepsilon$. But for this we need only to show that for a sufficiently small $\varepsilon > 0$, and uniformly in $n$,

$$\sup_{\varphi \in \Phi_\varepsilon} \left| \int_{1-\varepsilon}^1 \varphi(t) \gamma^T(t) \Gamma_t^{-1} \int_t^1 \gamma(s) \, du_n(s) \, dt \right|,$$

is arbitrarily small in probability. Denote the envelope function for $\varphi \in \Phi_\varepsilon$ by $\Psi$. Then, the above expression is bounded above by

$$\int_{1-\varepsilon}^1 |\Psi(t)| \left| \gamma^T(t) \Gamma_t^{-1} \int_t^1 \gamma(s) \, du_n(s) \right| dt.$$

However, bearing in mind that

$$E \left| \gamma^T(t) \Gamma_t^{-1} \int_t^1 \gamma(s) \, du_n(s) \right|^2 \leq \gamma^T(t) \Gamma_t^{-1} \gamma(t) \qquad \forall t \in [0,1],$$

we obtain that

$$E \int_{1-\varepsilon}^1 |\Psi(t)| \left| \gamma^T(t) \Gamma_t^{-1} \int_t^1 \gamma(s) \, du_n(s) \right| dt$$

$$= \int_{1-\varepsilon}^1 |\Psi(t)| E \left| \gamma^T(t) \Gamma_t^{-1} \int_t^1 \gamma(s) \, du_n(s) \right| dt$$

$$\leq \int_{1-\varepsilon}^1 |\Psi(t)| |\gamma^T(t) \Gamma_t^{-1} \gamma(t)|^{1/2} \, dt \leq \int_{1-\varepsilon}^1 \frac{1}{(1-t)^{1/2+\alpha}} \, dt.$$

The last integral can be made arbitrarily small for sufficiently small $\varepsilon$. $\quad\square$

Consequently, we obtain the following limit theorem for $b_n$. Recall, say from van der Vaart and Wellner (1996), that the family of Gaussian random variables $b(\varphi), \varphi \in L_2[0,1]$ with covariance function $Eb(\varphi)b(\varphi') = \int_0^1 \varphi(t)\varphi'(t) \, dt$ is called (function parametric) standard Brownian motion on $\Phi$ if $b(\varphi)$ is continuous on $\Phi$.

THEOREM 2.2. (i) *Let $\Phi$ be a Donsker class, that is, let $u_n \to_d u$ in $l_\infty(\Phi)$. Then, for every $\varepsilon > 0$,*

$$b_n \to_d b \qquad in \ l_\infty(\Phi \cap \Phi_\varepsilon),$$

*where $\{b(\varphi), \varphi \in \Phi\}$ is standard Brownian motion.*

(ii) *If the envelope function $\Psi(t)$ of (2.10) tends to positive (finite or infinite) limit at $t = 1$, then for the process (2.9) we have*

$$b_n \to_d b \qquad on \ [0, 1].$$

EXAMPLES. Here, we discuss some examples analyzing the behavior of the upper bound of (2.10) in the right tail. In all these examples we will see that not only the class of indicator functions satisfy (2.10) but also a class of unbounded functions $\varphi$ with $\varphi(s) = O((1-s)^{-\alpha}), \alpha < 1/2$, as $s \to 1$, satisfy this condition.

Consider logistic d.f. $F$ with the scale parameter 1, or equivalently $\psi_f(x) = 2F(x) - 1$. Then $h(x) = (1, 2F(x) - 1)^T$ or $\gamma(s) = (1, 2s - 1)^T$ and

$$\Gamma_s = (1-s) \begin{pmatrix} 1 & s \\ s & (1-2s+4s^2)/3 \end{pmatrix}, \qquad \det(\Gamma_s) = \frac{(1-s)^4}{3},$$

$$\Gamma_s^{-1} = \frac{3}{(1-s)^3} \begin{pmatrix} (1-2s+4s^2)/3 & -s \\ -s & 1 \end{pmatrix},$$

so that indeed $\gamma^T(s) \Gamma_s^{-1} \gamma(s) = 4(1-s)^{-1}$, for all $0 \le s < 1$.

Next, suppose $F$ is standard normal d.f. Because here $\psi_f(x) = x$, one obtains $h(x) = (1, x)^T$ and $\sigma_f^2(x) = xf(x) + 1 - F(x)$. Let $\mu(x) = f(x)/(1-F(x))$. Then,

$$\Gamma_{F(x)} = (1 - F(x)) \begin{pmatrix} 1 & \mu(x) \\ \mu(x) & x\mu(x) + 1 \end{pmatrix},$$

$$\Gamma_{F(x)}^{-1} = \frac{1}{(1-F(x))} \frac{1}{(x\mu(x) + 1 - \mu^2(x))} \begin{pmatrix} x\mu(x) + 1 & -\mu(x) \\ -\mu(x) & 1 \end{pmatrix}.$$

Hence

$$h^T(x) \Gamma_{F(x)}^{-1} h(x) = \frac{1}{(1-F(x))} \frac{(1 - x\mu(x) + x^2)}{(x\mu(x) + 1 - \mu^2(x))}.$$

Using asymptotic expansion for the tail of the normal d.f. [see, e.g., Feller (1957), page 179], for $\mu(x)$ we obtain

$$\mu(x) = \frac{x}{1 - S(x)} \qquad \text{where } S(x) = \sum_{i=1}^\infty \frac{(-1)^{i-1}(2i-1)!!}{x^{2i}} = \frac{1}{x^2} - \frac{3}{x^4} + \cdots.$$

From this one can derive that $(1 - x\mu(x) + x^2)/(x\mu(x) + 1 - \mu^2(x)) \sim 2$, $x \to \infty$, and therefore $h^T(x) \Gamma_{F(x)}^{-1} h(x) \sim 2(1 - F(x))^{-1}$, $x \to \infty$, or equivalently,

$$\gamma^T(s) \Gamma_s^{-1} \gamma(s) \sim 2(1-s)^{-1}, \qquad s \to 1.$$

Next, consider student $t_k$-distribution with fixed number of degrees of freedom $k$. In this case,

$$f(x) = \frac{1}{\sqrt{\pi k}} \frac{\Gamma((k+1)/2)}{\Gamma(k/2)} \frac{1}{(1+(x^2/k))^{(k+1)/2}},$$

$$\psi_f(x) = \frac{k+1}{k} \frac{x}{1+(x^2/k)}, \qquad x \in \mathbb{R}.$$

Using asymptotics for $k$ fixed and $x \to \infty$ we obtain [cf., e.g., Soms (1976)]

$$1 - F(x) \sim \frac{1+(x^2/k)}{x} f(x) \sim \frac{d_k}{k} \frac{1}{x^k}, \qquad d_k = \frac{1}{\sqrt{\pi}} \frac{\Gamma((k+1)/2)}{\Gamma(k/2)} k^{k/2}$$

$$f(x) \sim \frac{d_k}{x^{k+1}}, \qquad \psi_f(x) \sim \frac{(k+1)}{x}.$$

Consequently,

$$\Gamma_{F(x)} \sim \frac{d_k}{x^{k+2}} \begin{pmatrix} x^2/k & x \\ x & (k+1)^2/(k+2) \end{pmatrix},$$

$$\Gamma_{F(x)}^{-1} \sim \frac{x^k}{d_k} k(k+2) \begin{pmatrix} (k+1)^2/(k+2) & -x \\ -x & x^2/k \end{pmatrix},$$

$$h^T(x)\Gamma_{F(x)}^{-1}h(x) \sim \frac{2(k+1)}{d_k} x^k \sim \frac{2(k+1)}{k}[1-F(x)]^{-1}, \qquad x \to \infty,$$

or $\gamma^T(s)\Gamma_s^{-1}\gamma(s) \sim [2(k+1)/k](1-s)^{-1}$, as $s \to 1$.

The two values of $k = 1$ and $k = 2$ deserve special attention because mean and variance do not exist in these two cases. For $k = 1$, one obtains standard Cauchy distribution and, as seen above, the transformation *per ce* remains technically sound and the proposed test to fit the standard Cauchy distribution is valid as long as $m(x)$ is interpreted as some other conditional location parameter of $Y$, given $X = x$, such as conditional median, and as long as one has an estimator of this $m(x)$ satisfying (1.1). A similar comment applies when $k = 2$.

Finally, let $F$ be double exponential, or Laplace, d.f. with the density $f(x) = \alpha e^{-\alpha|x|}, \alpha > 0$. For $x > 0$ we get $h(x) = (1, \alpha)^T$ and $\gamma(s) = (1, \alpha)^T$, and $\Gamma_s$ becomes degenerate, equal to (2.6). Therefore again, see (2.7) with vector $b = \gamma(t)$, for $s > 1/2$, $\gamma^T(s)\Gamma_s^{-1}\gamma(s) = (1-s)^{-1}$.

Next, in this section we wish to clarify the question of a.s. continuity of $K_n$ and $K$ as linear functionals and thus justify the presence of tail condition (2.10). For this purpose it is sufficient to consider particular case, when $\gamma(s) = 1$ is one-dimensional and $\Gamma_s = 1 - s$. In this case

$$K_n(\varphi) = -\int_0^1 \varphi(s) \frac{u_n(s)}{1-s} \, ds, \qquad K(\varphi) = -\int_0^1 \varphi(s) \frac{u(s)}{1-s} \, ds.$$

The proposition below is of independent interest.

PROPOSITION 2.1.    (i) $K_n(\varphi)$ is continuous linear functional in $\varphi$ on $L_2[0,1]$ for every finite $n$.

(ii) However, the integral $\int_0^1 u^2(s)/(1-s)^2 \, ds$ is almost surely infinite. Moreover,

$$\frac{1}{-\ln(1-s)} \int_0^s \frac{u^2(t)}{(1-t)^2} \, dt \to_p 1 \qquad as \ s \to 1.$$

Therefore, $K(\varphi)$ is not continuous on $L_2[0,1]$.

REMARK 2.2.    It is easy to see that $E \int_0^1 u^2(s)/(1-s)^2 \, ds = \infty$, but this would not resolve the question of a.s. behavior of the integral and, hence, of $K$.

PROOF OF PROPOSITION 2.1.    (i) From the Cauchy–Schwarz inequality we obtain

$$|K_n(\varphi)| \le \left( \int_0^1 \varphi^2(s) \, ds \right)^{1/2} \left( \int_0^1 \frac{u_n^2(s)}{(1-s)^2} \, ds \right)^{1/2}$$

and the question reduces to whether the integral $\int_0^1 [u_n(s)/(1-s)]^2 \, ds$ is a.s. finite or not. However, it is, as even $\sup_s |u_n(s)/(1-s)|$ is a proper random variable for any finite $n$.

(ii) Recall that $u(s)/(1-s)$ is a Brownian motion: if $b$ denotes standard Brownian motion on $[0,\infty)$, then, in distribution,

$$\frac{u(t)}{1-t} = b\left( \frac{t}{1-t} \right) \qquad \forall t \in [0,1].$$

Hence, in distribution,

$$\int_0^s \frac{u^2(t)}{(1-t)^2} \, dt = \int_0^s b^2\left( \frac{t}{1-t} \right) dt = \int_0^\tau \frac{b^2(z)}{(1+z)^2} \, dz, \qquad \tau = s/(1-s).$$

Integrating the last integral by parts yields

$$\begin{aligned}
(2.11) \qquad \int_0^\tau \frac{b^2(z)}{(1+z)^2} \, dz &= -\frac{b^2(\tau)}{1+\tau} + 2 \int_0^\tau \frac{b(z)}{1+z} \, db(z) + \int_0^\tau \frac{1}{1+z} \, dz \\
&= -\frac{b^2(\tau)}{1+\tau} + 2 \int_0^\tau \frac{b(z)}{1+z} \, db(z) + \ln(1+z).
\end{aligned}$$

Consider the martingale

$$M(t) = \int_0^t \frac{b(z)}{1+z} \, db(z), \qquad t \ge 0.$$

Its quadratic variation process is

$$\langle M \rangle_t = \int_0^t \frac{b^2(z)}{(1+z)^2} \, dz.$$

Note that $\langle M \rangle_\tau$ equals the term on the left-hand side of (2.11). Divide (2.11) by $\ln(1+\tau)$ to obtain

$$\frac{\langle M \rangle_\tau}{\ln(1+\tau)} = -\frac{b^2(\tau)}{(1+\tau)\ln(1+\tau)} + 2\frac{M(\tau)}{\ln(1+\tau)} + 1.$$

The equalities

$$EM^2(t) = E\langle M \rangle_t = \int_0^t \frac{z}{(1+z)^2} \, dz = \ln(1+t) - \frac{1}{1+t}, \qquad Eb^2(t) = t,$$

imply that

$$\frac{b^2(\tau)}{(1+\tau)\ln(1+\tau)} = o_p(1) \quad \text{and} \quad \frac{M(\tau)}{\ln(1+\tau)} = o_p(1) \qquad \text{as } \tau \to \infty.$$

Hence, $\langle M \rangle_\tau / \ln(1+\tau) \to_p 1$, as $\tau \to \infty$.   □

**3. Power.** Consider, for the sake of comparison, the problem of fitting a distribution in the one sample location model up to an unknown location parameter. More precisely, consider the problem of testing that $X_1, \ldots, X_n$ is a random sample from $F(\cdot - \theta)$, for some $\theta \in \mathbb{R}$, against the class of all contiguous alternatives, that is, sequences of alternative distributions $A_n(\cdot - \theta)$ satisfying

$$\left( \frac{dA_n(x)}{dF(x)} \right)^{1/2} = 1 + \frac{1}{2\sqrt{n}} g(x) + r_n(x),$$

$$\int g^2(x) \, dF(x) < \infty, \qquad \int r_n^2(x) \, dF(x) = o\left( \frac{1}{n} \right).$$

As is known, and as can intuitively be understood, one should be interested only in the class of functions $g \in L_2(F)$ that are orthogonal to $\psi_f$:

$$(3.1) \qquad \int g(x)\psi_f(x) \, dF(x) = 0.$$

Indeed, as $g$ describes a functional "direction" in which the alternative $A_n$ deviates from $F$, if it has a component collinear with $\psi_f$,

$$g(x) = g_\perp(x) + c\psi_f(x), \qquad \int g_\perp(x)\psi_f(x) \, dF(x) = 0,$$

then infinitesimal changes in the direction $c\psi_f$ will be explained by, or attributed to, the infinitesimal changes in the value of parameter, that is,

"within" parametric family. Hence it cannot (and should not) be detected by a test for our parametric hypothesis. So, we assume that $g$ and $\psi_f$ are orthogonal, that is, (3.1).

Since $\theta$ remains unspecified, we still need to estimate it. Suppose $\bar{\theta}$ is its MLE under $F$ and consider empirical process $\bar{v}_n$ based on $\bar{e}_i = X_i - \bar{\theta}, i = 1, 2, \ldots, n$:

$$\bar{v}_n(x) = \sqrt{n}[\bar{F}_n(x) - F(x)], \qquad \bar{F}_n(x) = \frac{1}{n}\sum_{i=1}^n \mathbb{I}_{\{\bar{e}_i \leq x\}}.$$

One uses the empirical process $v_n$ in the case one assumes $\theta$ is known.

It is known [see, e.g., Khmaladze (1979)] that the asymptotic shift of $\bar{v}_n$ and $v_n$ under the sequence of alternatives $A_n$ with orthogonality condition (3.1) is the same and equals the function

$$G(x) = \int_{-\infty}^x g(y)\, dF(y).$$

However, the process $\bar{v}_n$ has uniform asymptotic representation

$$\bar{v}_n(x) = v_n(x) + f(x)\int \psi_f(y)\, dv_n(y) + o_p(1)$$

and, the main part on the right is orthogonal projection of $v_n$—see Khmaladze (1979) for a precise statement; see also Tjurin (1970). Heuristically speaking, it implies that the process $\bar{v}_n$ is "smaller" than $v_n$. In particular, variance of $\bar{v}_n(x)$ is bounded above by the variance of $v_n(x)$, for all $x$. Therefore, tests based on omnibus statistics, which typically measure an "overall" deviation of an empirical distribution function from $F$, or of empirical process from 0, will have better power if based on $\bar{v}_n$ than $v_n$. From a certain point of view this may seem a paradox, as it implies that, even if we know the parameter $\theta$, it would still be better to replace it by an estimator, because the power of many goodness of fit tests will thus increase. However, note that the integral in the last display has the same asymptotic distribution under hypothetical $F$ and alternatives $A_n$, and therefore the $v_n$ is "bigger" than $\bar{v}_n$ by the term, which is not useful in our testing problem.

Transformation of the process $\bar{v}_n$ asymptotically coincides with the process $w_n$ we study here, and moreover, the relationship between the two processes is one-to-one. Therefore, any statistic based on either one of these two processes will yield the same large sample inference.

With the process $\hat{v}_n$ the situation is the following: although it can be shown that the shift of this process under alternatives $A_n$ with orthogonality condition (3.1) is again function $G$, with general estimator $\hat{m}_n$ and, therefore, the general form of $R_n$, this process is not a transformation of $v_n$ only, and therefore is not its projection. In other words, it is not as "concentrated" as $\bar{v}_n$. The bias part of $R_n$ brings in additional randomization, not useful for the testing problem at hand. As a result, one will have less power in tests based on omnibus statistics from $\hat{v}_n$.

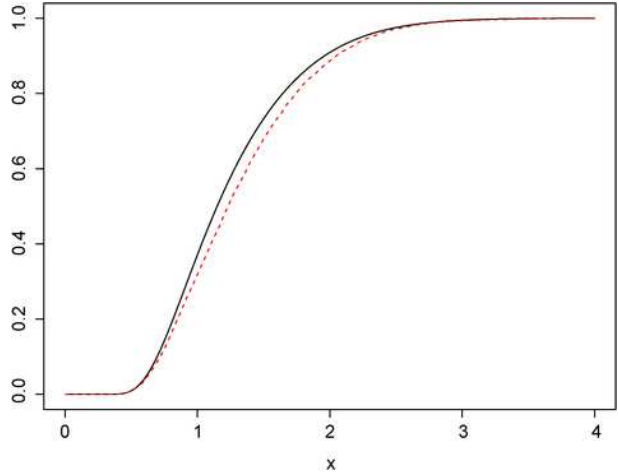

Fig. 1. *Null empirical d.f. (red dashed curve) and null limit d.f. (black curve) of $W_n$.*

We illustrate this by a simulation study. In this study we chose the regression model $Y = m(X) + e$, with $m(x) = e^x$, and covariate $X$ to be uniformly distributed on $[0,2]$. Let $F_0(\Psi)$ denote d.f. of a standardized normal (standardized double exponential) r.v. and $f_0(\psi)$ denote their densities. The problem is to test $H_0: F = F_0$, versus the alternatives $H_1: F \neq F_0$. In simulation below we chose a particular member of this alternative: $F_1 = 0.8F_0 + 0.2\Psi$. To estimate $m$, we used naive Nadaraya–Watson estimator

$$\hat{m}_n(x) = \sum_{i=1}^{n} Y_i \mathbb{I}_{\{X_i \in [x-a,x+a]\}} \Big/ \sum_{i=1}^{n} \mathbb{I}_{\{X_i \in [x-a,x+a]\}},$$

with $a = 0.04$. We shall compare the two tests based on $\hat{V}_n = \sup_x |\hat{v}_n(x)|$ and $W_n = \sup_x |w_n(x)|$. In all simulations, $n = 200$, repeated 10,000 times.

First, we generated null empirical d.f.'s of both statistics under the above set up. As seen in Figure 1, although the sample size $n = 200$ is not too big, the empirical null d.f. of $W_n$ is quite close to the d.f. of $\sup_x |b(F_0(x))|$, its limiting distribution. Empirical null d.f. of $\hat{V}_n$ is given in Figure 3.

To compare power of these tests, we generated 160 errors from $F_0$ and 40 from $\Psi$ and used the above set up to compute $\hat{V}_n$ and $W_n$. Figure 2 shows the hypothetical normal density $f_0$ versus the alternative mixture density $f_1 = 0.8f_0 + 0.2\psi$. Figure 3 describes empirical d.f.'s of $\hat{V}_n$ under $F_0$ and $F_1$ while Figure 4 gives the same entities for $W_n$.

Clearly, the alternative we consider, given that the sample size is only $n = 200$, should indeed be not easy to detect, especially by a test. Besides, as the difference between $F_0$ and $F_1$ occurs in the "middle" of the d.f. $F_0$, the alternative $F_1$ is of a nature, favorable for application of Komogorov–Smirnov test based on $\hat{v}_n$. However, Figures 3 and 4 show the effect we expected: distribution of $\hat{V}_n$ reacts to the alternative, that is, to the presence of double-exponential errors less than the distribution of $W_n$.

The above figures were computed with the window width $a = 0.04$. To assess the effect of window width on empirical power of these tests, we computed empirical power for additional values of $a = 0.08, 0.12$, at some empirical levels $\alpha$. Table 1 presents these numerical power values. In all cases one sees the empirical power of $W_n$ test to be larger than that of $\hat{V}_n$ test at all chosen levels $\alpha$, although for $a = 0.04$, this difference is far more significant than in the other two cases. Critical values used in this comparison were estimated from their respective empirical null distributions. These are not isolated findings—more examples can be found in Brownrigg (2008).

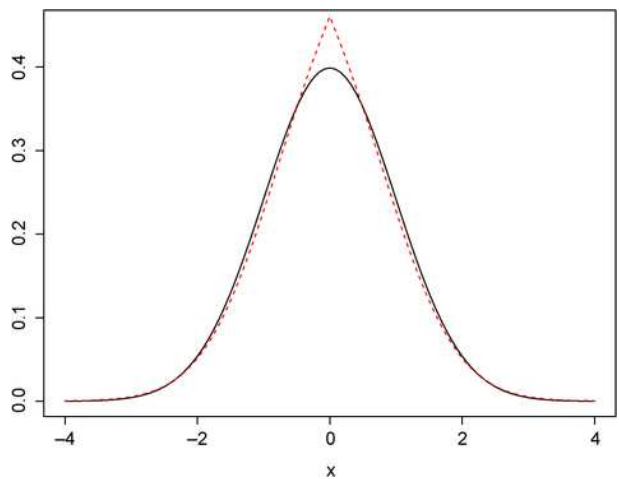

Fig. 2. $f_0$ (dark curve) and $0.8f_0 + 0.2\psi$ (red dashed curve).

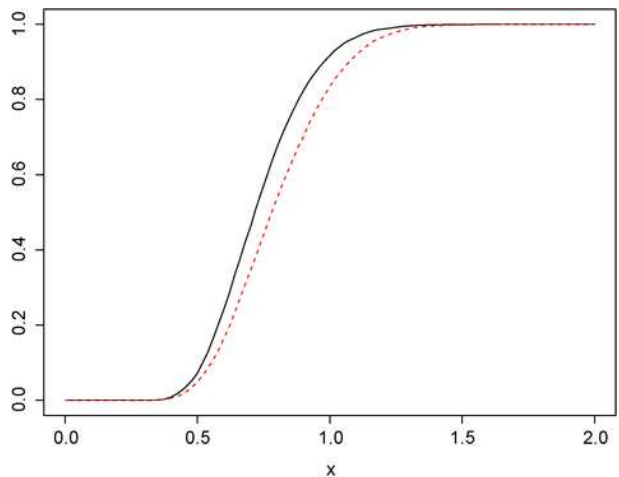

Fig. 3. Empirical d.f.'s of $\hat{V}_n$ under $H_0$ (black curve) and $H_1$ (red dashed curve).

TABLE 1
*Empirical power of $\hat{V}_n$ and $W_n$ tests*

| $f$ and $a$ | $\alpha$ | $\hat{V}_n$ | $W_n$ |
|---|---|---|---|
| $f_1$, $a = 0.04$ | 0.10 | 0.1904 | 0.3168 |
| | 0.05 | 0.1154 | 0.1920 |
| | 0.025 | 0.0625 | 0.1114 |
| | 0.01 | 0.0260 | 0.0523 |
| $f_1$, $a = 0.08$ | 0.10 | 0.1838 | 0.2115 |
| | 0.05 | 0.1081 | 0.1242 |
| | 0.025 | 0.0680 | 0.0744 |
| | 0.01 | 0.0325 | 0.0450 |
| $f_1$, $a = 0.12$ | 0.10 | 0.1837 | 0.1960 |
| | 0.05 | 0.1085 | 0.1150 |
| | 0.025 | 0.0619 | 0.0760 |
| | 0.01 | 0.0301 | 0.0480 |

Returning to general discussion on power, we must add that with the estimator $\hat{m}_n$ used by Müller, Schick and Wefelmeyer, and therefore, with their simple form of $R_n$, the process $\hat{v}_n$ is again asymptotically a projection, although in general a skew one, of the process $v_n$. As described in Khmaladze (1979), it is asymptotically in one-to-one relationship with the process $\bar{v}_n$, and, therefore $w_n$. Hence, the large sample inference drawn from a statistic based on $\hat{v}_n$ is, in this case, also equivalent to that drawn from the analogous statistic based on either of the other two, and the only difference between this processes is that $\hat{v}_n$ and $\bar{v}_n$ are not asymptotically distribution free, while $w_n$ is.

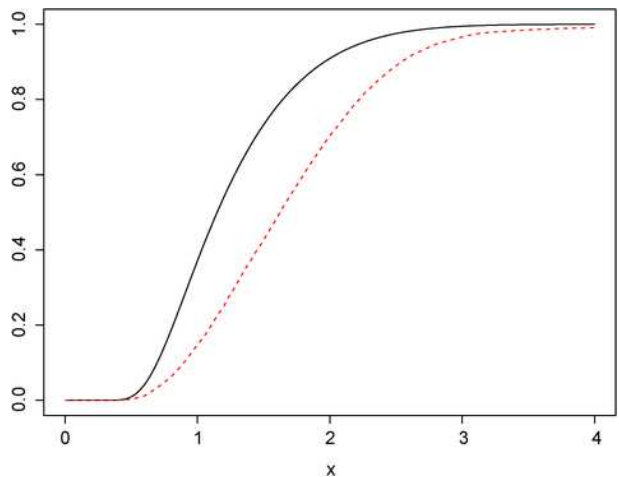

FIG. 4. *Empirical d.f.'s of $W_n$ under $H_0$ (black curve) and $H_1$ (red dashed curve).*

**4. Weak convergence of $w_n$.** In this section we prove weak convergence for the process $w_n$, given by (2.2) and (2.3). In view of (2.5), (2.9) and the fact that the weak convergence of the first part in the right-hand side of (2.5) was proved in Theorem 2.1, it suffices to show that the process $\eta_n$ of (2.5) is asymptotically small. Being the transformation of "small" process $\xi_n$, the smallness of $\eta_n$ is plausible. However, the transformation $K(\cdot, \xi_n)$ is not continuous in $\xi_n$ in uniform metric. Indeed, although in the integration by parts formula

$$\int_t^1 \gamma(s)\, d\xi_n(F^{-1}(s)) = \xi_n(F^{-1}(s))\gamma(s)|_{s=t}^1 - \int_t^1 \xi_n(F^{-1}(s))\, d\gamma(s),$$

we can show, that $\xi_n(F^{-1}(1))\gamma(1) = 0$, the integral on the right-hand side is not necessarily small if $\gamma(t)$ is not bounded at $t = 1$, as happens to be true for normal d.f. $F$ where the second coordinate of $\gamma(t)$ is $F^{-1}(t)$. Therefore, one cannot prove the smallness of $\eta_n$ in sufficient generality, using only uniform smallness of $\xi_n$.

If we use, however, quite mild additional assumption on the right tail of $\xi_n$, or rather of $\hat{v}_n$ and $f$, we can obtain the weak convergence of $w_n$ basically under the same conditions as in Theorem 2.2. Namely, assume that for some positive $\beta < 1/2$,

$$(4.1) \qquad \sup_{y>x} \frac{|\hat{v}_n(y)|}{(1-F(y))^\beta} = o_p(1) \qquad \text{as } x \to \infty,$$

uniformly in $n$. Note that the same condition for $v_n$ is satisfied for all $\beta < 1/2$.

Denote tail expected value and variance of $\psi_f(e_1)$ by

$$E[\psi_f|x] = E[\psi_f(e_1)|e_1 > x], \qquad \text{Var}[\psi_f|x] = \text{Var}[\psi_f(e_1)|e_1 > x].$$

Now we formulate two more conditions on $F$.

(a) For any $\varepsilon > 0$ the function $\psi_f(F^{-1})$ is of bounded variation on $[\varepsilon, 1-\varepsilon]$ and for some $\varepsilon > 0$ it is monotone on $[1-\varepsilon, 1]$.

(b) For some $\delta > 0$, $\varepsilon > 0$ and some $C < \infty$,

$$\frac{(\psi_f(x) - E[\psi_f|x])^2}{\text{Var}[\psi_f|x]} < C(1-F(x))^{-2\delta} \qquad \forall x: F(x) > 1-\varepsilon.$$

Note that in terms of the above notation,

$$(4.2) \quad \gamma(t)^T \Gamma_t^{-1} \gamma(t) = \frac{1}{1-F(x)}\left[1 + \frac{(\psi_f(x) - E[\psi_f|x])^2}{\text{Var}[\psi_f|x]}\right], \qquad t = F(x).$$

Hence, condition (b) is equivalent to

$$(4.3) \qquad \gamma(t)^T \Gamma_t^{-1} \gamma(t) \le C(1-t)^{-1-2\delta} \qquad \forall t > 1 - \varepsilon.$$

Condition (b) is easily satisfied in all examples of Section 2, even with $\delta = 0$.

Our last condition is as follows.

(c) For some $C < \infty$ and $\beta > 0$ as in (4.1),

$$\left| \int_x^\infty [1 - F(y)]^\beta \, d\psi_f(y) \right| \le C |\psi_f(x) - E[\psi_f | x]|.$$

Condition (c) is also easily satisfied in all examples of Section 2, even for arbitrarily small $\beta$.

For example, for logistic distribution, with $t = F(x)$, $\psi_f(x) = 2t - 1$ and

$$\left| \int_x^\infty [1 - F(y)]^\beta \, d\psi_f(y) \right| = 2 \int_t^1 (1-s)^\beta \, ds = \frac{2}{\beta+1}(1-t)^{\beta+1},$$

while $|\psi_f(x) - E[\psi_f | x]| = (1 - t)$ and their ratio tends to 0, as $t \to 1$. For normal distribution,

$$\int_x^\infty [1 - F(y)]^\beta \, d\psi_f(y) \sim \int_x^\infty \frac{1}{y^\beta} f^\beta(y) \, dy \le \frac{1}{x} \int_x^\infty y^{1-\beta} f^\beta(y) \, dy,$$

while

$$|\psi_f(x) - E[\psi_f | x]| = \left| x - \frac{f(x)}{1 - F(x)} \right| \sim \frac{x}{x^2 - 1}, \qquad x \to \infty,$$

and the ratio again tends to 0, as $x \to \infty$.

Recall the notation

$$K(\varphi, \xi_n) = \int_0^1 \varphi(t) \gamma(t)^T \Gamma_t^{-1} \int_t^1 \gamma(s) \xi_n(F^{-1}(ds)) \, dt$$

and for a given indexing class $\Phi$ of functions from $L_2[0,1]$ let $\Phi \circ F = \{\varphi(F(\cdot)), \varphi \in \Phi\}$.

THEOREM 4.1. (i) *Suppose conditions (4.1) and* (a)–(c) *are satisfied with* $\beta > \delta$. *Then, on the class* $\Phi_\varepsilon$ *as in Theorem 2.1 but with* $\alpha < \beta - \delta$, *we have*

$$\sup_{\varphi \in \Phi_\varepsilon} |K(\varphi, \xi_n)| = o_p(1), \qquad n \to \infty.$$

*Therefore, if* $\Phi$ *is a Donsker class, then, for every* $\varepsilon > 0$,

$$w_n \to_d b \qquad in \ l_\infty(\Phi \cap \Phi_\varepsilon \circ F),$$

*where* $\{b(\varphi), \varphi \in \Phi\}$ *is standard Brownian motion.*

(ii) *If, in addition,* $\delta \le \alpha$, *then for the time transformed process* $w_n(F^{-1}(\cdot))$ *of (2.2), we have*

$$w_n(F^{-1}(\cdot)) \to_d b(\cdot) \qquad in \ D[0,1].$$

PROOF.   Note, that

$$\gamma(t)^T \Gamma_t^{-1}(0, a)^T = \frac{1}{1 - F(x)} \frac{(\psi_f(x) - E[\psi_f | x]) a}{\mathrm{Var}[\psi_f | x]}, \qquad t = F(x), \ \forall a \in \mathbb{R}.$$

Use this equality for $a = \int_t^1 (1 - s)^\beta \, d\psi_f(F^{-1}(s))$. Then condition (c) implies that

$$(4.4) \qquad |\gamma(t)^T \Gamma_t^{-1}(0, a)^T| \le C \gamma(t)^T \Gamma_t^{-1} \gamma(t) \qquad \forall t < 1.$$

Now we prove the first claim.

(i) Use the notation $\xi_n'(t) = \xi_n(x)$ with $t = F(x)$. Since we expect singularities at $t = 0$ and, especially, at $t = 1$ in both integrals in $K(\varphi, \xi_n)$ we will isolate the neighborhood of these points and consider it separately. Mostly we will take care of the neighborhood of $t = 1$. The neighborhood of $t = 0$ can be treated more easily (see below). First assume $\Gamma_t^{-1}$ nondegenerate for all $t < 1$. Then,

$$
\begin{aligned}
(4.5) \qquad & \int_0^1 \varphi(t) \gamma(t)^T \Gamma_t^{-1} \int_t^1 \gamma(s) \xi_n'(ds) \, dt \\
& = \int_0^{1-\varepsilon} \varphi(t) \gamma(t)^T \Gamma_t^{-1} \int_t^{1-\varepsilon} \gamma(s) \xi_n'(ds) \, dt \\
& \quad + \int_0^{1-\varepsilon} \varphi(t) \gamma(t)^T \Gamma_t^{-1} \int_{1-\varepsilon}^1 \gamma(s) \xi_n'(ds) \, dt \\
& \quad + \int_{1-\varepsilon}^1 \varphi(t) \gamma(t)^T \Gamma_t^{-1} \int_t^1 \gamma(s) \xi_n'(ds) \, dt.
\end{aligned}
$$

Consider the third summand on the right-hand side. First note that, when proving that it is small, we can replace $\xi_n$ by the difference $\hat{v}_n - v_n$ only. Indeed, since $df(F^{-1}(s)) = \psi_f(x) f(x) \, dx$, according to (2.4) the integral

$$\int_{1-\varepsilon}^1 \varphi(t) \gamma(t)^T \Gamma_t^{-1} \int_t^1 \gamma(s) \, df(F^{-1}(s)) \, dt$$

is the second coordinate of $\int_{1-\varepsilon}^1 \varphi(t) \gamma(t) \, dt$, and is small for $\varepsilon$ small anyway. Monotonicity of $\psi_f(F^{-1})$ guaranteed by assumption (a) and (2.1) justify integration by parts of the inner integral in the following derivation.

$$
\begin{aligned}
& \int_{1-\varepsilon}^1 \varphi(t) \gamma(t)^T \Gamma_t^{-1} \int_t^1 \gamma(s) \hat{u}_n(ds) \, dt \\
& = \int_{1-\varepsilon}^1 \varphi(t) \gamma(t)^T \Gamma_t^{-1} \left[ -\gamma(t) \hat{u}_n(t) - \int_t^1 \hat{u}_n(s) \, d\gamma(s) \right] dt.
\end{aligned}
$$

Assumption (2.10) on $\varphi$ and (4.3) imply

$$\left| \int_{1-\varepsilon}^1 \varphi(t)\gamma(t)^T \Gamma_t^{-1} \gamma(t) \hat{u}_n(t)\, dt \right|$$

$$\leq C \int_{1-\varepsilon}^1 [\gamma(t)^T \Gamma_t^{-1} \gamma(t)]^{1/2} \frac{1}{(1-t)^{1/2+\alpha-\beta}}\, dt \sup_{t>1-\varepsilon} \frac{|\hat{u}_n(t)|}{(1-t)^\beta}$$

$$\leq C \int_{1-\varepsilon}^1 \frac{1}{(1-t)^{1+\alpha+\delta-\beta}}\, dt \sup_{t>1-\varepsilon} \frac{|\hat{u}_n(t)|}{(1-t)^\beta},$$

which is small for small $\varepsilon$ as soon as $\alpha < \beta - \delta$.

Now, note that $\int_t^1 \hat{u}_n(s)\, d\gamma(s) = (0, \int_t^1 \hat{u}_n(s)\, d\psi_f(F^{-1}(s)))^T$. Using monotonicity of $\psi_f(F^{-1})$ for small enough $\varepsilon$ we obtain, for all $t > 1 - \varepsilon$,

$$(4.6) \quad \left| \int_t^1 \hat{u}_n(s)\, d\psi_f(F^{-1}(s)) \right| < C \left| \int_t^1 (1-s)^\beta\, d\psi_f(F^{-1}(s)) \right| \sup_{s>1-\varepsilon} \frac{|\hat{u}_n(s)|}{(1-s)^\beta}.$$

Therefore, using (4.4), for the double integral we obtain

$$\left| \int_{1-\varepsilon}^1 \varphi(t)\gamma(t)^T \Gamma_t^{-1} \int_t^1 \hat{u}_n(s)\, d\gamma(s)\, dt \right|$$

$$\leq C \int_{1-\varepsilon}^1 |\varphi(t)| \gamma(t)^T \Gamma_t^{-1} \gamma(t)\, dt \sup_{s>1-\varepsilon} \frac{|\hat{u}_n(s)|}{(1-s)^\beta}$$

and the integral on the right-hand side, as we have seen above, is small as soon as $\alpha < \beta - \delta$. The same conclusion is true for $\hat{u}_n$ replaced by $u_n$.

Since (4.6) implies the smallness of

$$\int_{1-\varepsilon}^1 \hat{u}_n(s)\, d\psi_f(F^{-1}(s)) \quad \text{and} \quad \int_{1-\varepsilon}^1 u_n(s)\, d\psi_f(F^{-1}(s)),$$

to prove that the middle summand on the right-hand side of (4.5) is small one needs only finiteness of $\psi_f(x)$ in each $x$ with $0 < F(x) < 1$, which follows from (a). This and uniform in $x$ smallness of $\xi_n$ proves smallness of the first summand as well.

The smallness of integrals

$$\int_0^\varepsilon \varphi(t)\gamma(t)^T \Gamma_t^{-1} \gamma(t) \int_t^1 \gamma(s)\xi_n'(ds)\, dt$$

follows from $\Gamma_t^{-1} \sim \Gamma_0^{-1}$ and square integrability of $\varphi$ and $\gamma$.

If $\Gamma_t^{-1}$ becomes degenerate after some $t_0$, for these $t$ we get

$$\gamma(t)^T \Gamma_t^{-1} \int_t^1 \gamma(s)\xi_n'(ds) = \frac{\xi_n'(t)}{1-t}$$

and the smallness of all tail integrals easily follows for our choice of the indexing functions $\varphi$.

(ii) Since for $\delta \leq \alpha$ the envelope function $\Psi(t)$ of (2.10) satisfies inequality

$$\Psi(t) \geq (1-t)^{\delta-\alpha},$$

it has positive finite or infinite lower limit at $t = 1$. But then it is possible to choose as an indexing class the class of indicator functions $\varphi(t) = \mathbb{I}_{\{t \leq \tau\}}$ and the claim follows. □

REMARK 4.1 (Computational formula). We present here a computational formula for $w_n$. Let $\mathcal{G}(x) = \int_{y \leq x} \Gamma_{F(y)}^{-1} h(y) \, dF(y)$. Then, using (2.3) and (2.4) one obtains

$$w_n(x) = n^{-1/2} \sum_{i=1}^{n} [I(\hat{e}_i \leq x) - h(\hat{e}_i)^T \mathcal{G}(x \wedge \hat{e}_i)], \qquad x \in \mathbb{R}.$$

Thus to implement test based on $\sup_x |w_n(x)|$, one needs to evaluate $\mathcal{G}$ and compute $\max_{1 \leq j \leq n} |w_n(\hat{e}_{(j)})|$, where $\hat{e}_{(j)}$, $1 \leq j \leq n$, are the order statistics of $\hat{e}_j, 1 \leq j \leq n$.

REMARK 4.2 (Testing with an unknown scale). Here, we shall describe an analog of the above transformation suitable for testing the hypothesis $H_{sc}$ that the common d.f. of the errors $e_i$ is $F(x/\sigma)$, for all $x \in \mathbb{R}$, and for some $\sigma > 0$. Let $\phi_f(x) = 1 + x\psi_f(x)$ and $h_\sigma(x) = (1, \sigma^{-1}\psi_f(x), \sigma^{-1}\phi_f(x))^T$. Then analog of the vector $h(x)$ here is $h_\sigma(x/\sigma)$ and that of $\Gamma_t$ is

$$\Gamma_{t,\sigma} = \int_{y \geq x/\sigma} h_\sigma(y) h_\sigma(y)^T \, dF(y), \qquad t = F\left(\frac{x}{\sigma}\right).$$

This is the same matrix as given in Khmaladze and Koul (2004), page 1013. Akin to the function $K(x, \nu)$ define

$$K_\sigma(x, \nu) = \int_{-\infty}^{x/\sigma} h_\sigma^T(y) \Gamma_{F(y),\sigma}^{-1} \int_{y/\sigma}^{\infty} h_\sigma(z) \, d\nu(z\sigma) \, dF(y), \qquad x \in \mathbb{R}.$$

Analog of Lemma 2.1 continues to hold for each $\sigma > 0$, and hence this function is well defined for all $x \in \mathbb{R}, \sigma > 0$.

Let $\hat{\sigma}$ be a $n^{1/2}$-consistent estimator of $\sigma$ based on $\{(X_i, Y_i), 1 \leq i \leq n\}$. Let $\tilde{F}_n(x)$ be the empirical d.f. of the residuals $\tilde{e}_i = \hat{e}_i/\hat{\sigma}$ and let $\tilde{v}_n = n^{1/2}[\tilde{F}_n - F]$. Then the analog of $w_n$ suitable for testing $H_{sc}$ is

$$\tilde{w}_n(x) = n^{1/2}[\tilde{F}_n(x) - K_{\hat{\sigma}}(x, \tilde{F}_n)] = \tilde{v}_n(x) - K_{\hat{\sigma}}(x, \tilde{v}_n).$$

Under conditions analogous to those given in Section 4 above, one can verify that the conclusions of Theorem 4.1 continue to hold for $\tilde{w}_n$ also.

If we let $\mathcal{G}_\sigma(x) = \int_{y \leq x/\sigma} \Gamma_{F(y),\sigma}^{-1} h_\sigma(y) \, dF(y)$, then, one can rewrite

$$\tilde{w}_n(x) = n^{-1/2} \sum_{i=1}^{n} [I(\tilde{e}_i \leq x) - h_{\hat{\sigma}}(\tilde{e}_i)^T \mathcal{G}_{\hat{\sigma}}(x \wedge \tilde{e}_i)], \qquad x \in \mathbb{R}.$$

Hence, $\sup_x |\widetilde{w}_n(x)| = \max_{1 \le j \le n} |\widetilde{w}_n(\widetilde{e}_{(j)})|$, where $\widetilde{e}_{(j)}$, $1 \le j \le n$, are the order statistics of $\widetilde{e}_j$, $1 \le j \le n$.

**Acknowledgments.** The authors thank Dr. Ray Brownrigg for his invaluable help with the simulation studies appearing in Section 3, to a referee and the Associate Editor for their comments, and to the Co-Editor Bernard Silverman for his patience.

School of Mathematics, Statistics
  and Operations Research
Victoria University of Wellington
PO Box 600, Wellington
New Zealand
E-mail: Estate.Khmaladze@vuw.ac.nz

Department of Statistics
  and Probability
Michigan State University
East Lansing, Michigan 48824-1027
USA
E-mail: koul@stt.msu.edu